\documentclass{article}
\usepackage{spconf,amsmath,graphicx}
\usepackage{amssymb}  
\usepackage[utf8]{inputenc}
\usepackage[amsmath,thmmarks,hyperref]{ntheorem} 
\usepackage{bbm}               
\usepackage{mathtools}         
\usepackage{stmaryrd}          
\usepackage{paralist}          
\usepackage{aliascnt}          
\usepackage{subcaption} 
\usepackage{mathrsfs}

\usepackage{textcomp}

\usepackage{color}             
\graphicspath{{pics/}}

\newcommand\MTkillspecial[1]{
\bgroup
\catcode`\&=9
\let\\\relax%
\scantokens{#1}%
\egroup
}

\renewcommand{\mathbb}[1]{\mathbbm{#1}}

\theoremheaderfont{\normalfont\itshape}
\theorembodyfont{\normalfont}

\newaliascnt{proposition}{lemma}

\aliascntresetthe{proposition}

\newaliascnt{thm}{lemma}

\aliascntresetthe{thm}

\newaliascnt{corollary}{lemma}

\aliascntresetthe{corollary}

\newaliascnt{definition}{lemma}

\aliascntresetthe{definition}

\newaliascnt{claim}{lemma}

\aliascntresetthe{claim}

\theorembodyfont{\normalfont}

\newaliascnt{example}{lemma}

\aliascntresetthe{example}

\newaliascnt{remark}{lemma}

\aliascntresetthe{remark}

\newaliascnt{question}{lemma}

\aliascntresetthe{question}

\newaliascnt{conjecture}{lemma}

\aliascntresetthe{conjecture}

\theorembodyfont{\normalfont}
\theoremnumbering{Roman}

\def\theorem@checkbold{}

\allowdisplaybreaks

\let\originalleft\left
\let\originalright\right
\renewcommand{\left}{\mathopen{}\mathclose\bgroup\originalleft}
\renewcommand{\right}{\aftergroup\egroup\originalright}

\newcommand{\I}              {\mathrm{i}}
\newcommand{\E}              {\mathrm{e}}

\newcommand{\algebra}[1]      {\mathscr{#1}}

\newcommand{\Fourier}{\mathcal{F}}

\DeclarePairedDelimiter\ideal\langle\rangle
\reDeclarePairedDelimiterInnerWrapper\ideal{star}{
\mathopen{#1\vphantom{\MTkillspecial{#2}}\kern-\nulldelimiterspace\right.}
#2
\mathclose{\left.\kern-\nulldelimiterspace\vphantom{\MTkillspecial{#2}}#3}}

\newcommand{\DFT}{\mathsf{DFT}}

\makeatletter
\let\@@mod\mod
\DeclareRobustCommand{\mod}{\@ifstar\@mods\@@mod}
\def\@mods{\mkern4mu{\operator@font mod}\mkern6mu}
\makeatother

\title{The discrete cosine transform on triangles}
\twoauthors
  {Bastian Seifert\sthanks{This work is partially
    supported by the European Regional Development Fund (ERDF)}}
	{Ansbach University of Applied Sciences\\
	Faculty of Engineering Sciences\\
	Ansbach, Germany \\
      bastian.seifert@hs-ansbach.de}
  {Knut Hüper}
	{University of Würzburg\\
	Institute of Mathematics\\
	Würzburg, Germany\\
      hueper@mathematik.uni-wuerzburg.de}
\begin{document}
%
\maketitle
\begin{abstract}
    The discrete cosine transform is a valuable tool in analysis of
    data on undirected rectangular grids, like images. In this paper
    it is shown how one can define an analogue of the discrete cosine
    transform on triangles. This is done by combining algebraic signal
    processing theory with a specific kind of multivariate Chebyshev
    polynomials. Using a multivariate Christoffel-Darboux formula it
    is shown how to derive an orthogonal version of the transform.
\end{abstract}
\section{Introduction}
\label{sec:introduction}%

Triangulation of surfaces is a well-established method to discretize
surfaces and investigate geometric data. Using signal processing
techniques on surfaces requires in general either knowledge of the
surface itself or one must require that signals on the boundary of a
parametrisation match due to the periodic boundary conditions of the
discrete Fourier transform.

In this work we discuss another approach to a spectral transform on
triangles. This transform mimics the cosine transform on
one-dimensional data. Here the triangles are implemented as lattices
of triangles (not to be confused with triangular lattices as in
\cite{Bodner.Patera.Szajewska:2017a}). The usage of the presented
transform can be advantageous compared to the technique used
in~\cite{Ryland.Munthe-Kaas:2011} for spectral apprxoimation of
triangles as one does not need to handle a cusped region e.g. one
similar to the deltoid.

The derivation is based on algebraic signal processing (ASP)
theory~\cite{Pueschel.Moura:2008a,Pueschel.Moura:2008b}. In ASP one
reveals the algebraic principles underlying discrete signal processing
techniques.

These techniques were used in \cite{Seifert.Hueper.Uhl:2018a} to
derive a cosine transform together with its fast algorithm on the
face-centered cubic lattice and in \cite{Pueschel.Roetteler:2007} on
the hexagonal lattice. Both approaches relied on multivariate
Chebyshev polynomials. Multivariate Chebyshev polynomials are much
lesser known than their univariate counterparts, even though they
share many of their nice properties. In fact multivariate Chebyshev
polynomials can be deduced from principles well established in Lie
theory~\cite{Hoffman.Withers:1988}. Their construction is based on a
generalized cosine, which resembles the folding of some special
region. These regions are in one-to-one correspondence to finite
reflection groups and can be enumerated using the notion of
Coxeter-Dynkin diagrams~\cite{Conway.Sloane:1999}. The Coxeter-Dynkin
diagrams can be partitioned into series denoted by $A_n, B_n, C_n,$
and $D_n$ (and 5 additional special cases)~\cite{Conway.Sloane:1999}.
In~\cite{Seifert.Hueper.Uhl:2018a,Pueschel.Roetteler:2007} $A_n$-type
Chebyshev polynomials were used. In this work we rely on multivariate
Chebyshev polynomials of $B_2$-type. We will use an elementary
approach for their construction. Consequently no knowledge about Lie
theory is needed. After recalling the basic principles of ASP in
Sect.~\ref{sec:AlgebraicSignalProcessing} we will define the
$B_2$-Chebyshev polynomials, discuss some of their properties and
investigate the associated signal model and transform in
Sect.~\ref{sec:CosineOnTriangles}.

The Gauß-Jacobi procedure to derive unitary and orthogonal versions of
signal transforms~\cite{Yemini.Pearl:1979a} is connected to the
Christoffel-Darboux formula for orthogonal polynomials. Even though
there is a multivariate version available~\cite{Xu:1993a} we are not
aware of its usage in signal processing. In
Sect.~\ref{sec:ChristoffelDarbouxForInverse} we will use this
multivariate Christoffel-Darboux formula to derive an orthogonal
version of the transform defined in this paper. To our knowledge this
is the first time this formula is used in signal processing.

\section{Algebraic signal processing}
\label{sec:AlgebraicSignalProcessing}%

In this section the methods of algebraic signal processing theory are
illustrated on discrete signal processing with finite-time. This in
turn motivates the development of the tools in the next sections.

In finite time discrete signal processing a set of numbers
$s = (s_0, \dots, s_{n-1}) \in \mathbb{C}^n$ is called a signal if it
is periodically extended. That is one has $s_N = s_{N \mod* n}$ for any
$N \in \mathbb{Z}$. The finite $z$-transform associates to a signal a
polynomial in $x = z^{-1}$
\begin{equation}
    \label{eq:zTransform}
    (s_0,\dots,s_{n-1}) \mapsto \sum_{i=0}^{n-1} s_i x^i. 
\end{equation}
The periodic extension is captured by considering the polynomials
modulo $x^n - 1$, i.e.\, one requires $x^n = 1$. The set of
polynomials modulo $x^n - 1$ (or more precisely modulo the ideal
$\ideal{x^n - 1}$) is denoted by
$\mathbb{C}[x] \big/ \ideal{x^n - 1}$. So the finite $z$-transform is
a map
$\Phi \colon \mathbb{C}^n \longrightarrow \mathbb{C}[x] \big/
\ideal{x^n - 1}$.

A central concept in signal processing is the notion of shift. In
the $z$-domain the shift of finite time signal processing can be
realized as multiplication by $x$
\begin{equation}
    \label{eq:zDomainShift}
    x \cdot \Phi(s) = x \sum_{i=0}^{n-1} s_i x^i = \sum_{i=0}^{n-1}
    s_{i-1 \mod* n} x^i,
\end{equation}
and results in a delay of the signal.

Filters can be described in the $z$-domain as polynomials in the shift
$x$, i.e. a filter $h$ is of the form $h = \sum_{i=0}^{n-1} h_i x^i$.
Then filtering is just multiplication in
$\mathbb{C}[x] \big/ \ideal{x^n - 1}$
\begin{equation}
    \label{eq:zDomainFiltering}
    h \Phi(s) = \left( \sum_{i=0}^{n-1} h_i x^i \right) \left(
        \sum_{i=0}^{n-1} s_i x_i \right) \mod* x^n - 1.
\end{equation}
This notion of filtering is thus circular convolution. Furthermore one
can multiply filters modulo $x^n -1$ by one another to produce new
filters.

Mathematically this turns the set of filters into a polynomial
algebra. Since it is not meaningful to multiply signals we do not have
an algebra strucutre on the signals. But through the filtering
operation they get the structure of a module over the filter algebra.
The $z$-transform is a bijection between the data equipped with no
structure to a representation of the data equipped with algebraic
signal structure. So an algebraic signal model consists of a triple
$(\algebra{A},M,\Phi)$, where $\algebra{A}$ is an algebra, $M$ an
$\algebra{A}$-module and $\Phi \colon \mathbb{C}^n \to M$ a bijection.

From this data one obtains automatically a notion of Fourier
transform. The polynomial $x^n - 1 = \prod_k x - \E^{2 \pi \I k/n}$
decomposes into linear factors of its zeros. This leads, using the
Chinese remainder theorem, to a decomposition of the module into
irreducible submodules, i.e. one has an isomorphism
\begin{equation}
    \label{eq:FourierDecompositionModule}
    \mathbb{C}[x] \big/ \ideal{x^n - 1}
    \longrightarrow
    \bigoplus_k \mathbb{C}[x] \big/ \ideal{x - \E^{2 \pi \I k/n}}.
\end{equation}
Any matrix realizing this isomorphism is called a Fourier transform of
the model. For example if one chooses $(1,x,\dots,x^{n-1})$ as basis
for $\mathbb{C}[x] \big/ \ideal{x^n - 1}$ and $(1)$ in each irreducible
submodule one obtains the discrete Fourier transform matrix
$\DFT_n = \left[ \E^{2 \pi \I k j/n} \right]_{k,j}$.

\begin{figure}
    \centering
    \includegraphics[width=0.45\textwidth]{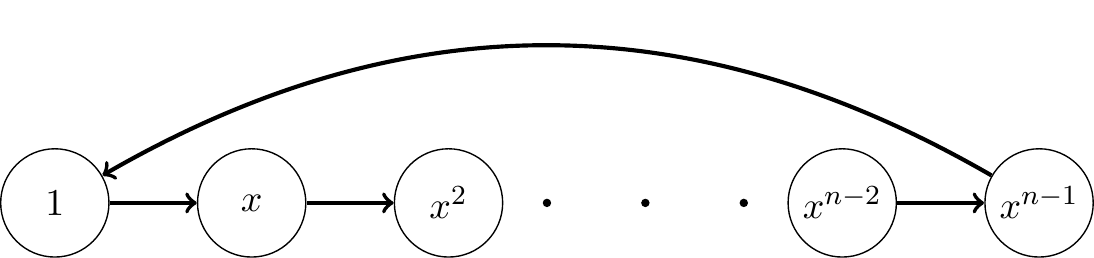}
    \caption{Visualization graph for the finite time discrete signal
      processing model. }
    \label{fig:DFTModelEngineer}
\end{figure}

A signal model can furthermore be visualized by a graph. This
visualization graph is obtained by adding a node for each basis
element of the module. Then an edge from one basis element to another
is added if multiplication by the generators of the algebra lead to a
linear combination containing the end point of the edge. For the
finite time discrete signal processing model this visualization is
shown in Fig.~\ref{fig:DFTModelEngineer}. It shows the periodic
extension of the underlying signal model. The fact that the graph is
directed gives rise to a time-model. If one replaces $x^n - 1$ for
example with the Chebyshev polynomial $T_n$ and the basis
$\{1,x,\dots,x^{n-1}\}$ by the basis $\{T_0,\dots,T_{n-1}\}$ one
obtains an undirected graph (and a discrete cosine transform), which
corresponds to a space model~\cite{Pueschel.Moura:2008c}.

\section{Cosine transform on triangles}
\label{sec:CosineOnTriangles}%

We now investigate a generalization of Chebyshev polynomials to
polynomials in two variables. The generalization is a special case,
$B_2$, of a whole family of generalized Chebyshev polynomials
associated to Lie theory. Here we use a down-to-earth approach to
these polynomials, for the general theory
see~\cite{Hoffman.Withers:1988}.

As in the univariate case, where one has
$T_n(x) = T_n(\cos \theta) = \cos n \theta$, it is much more
convenient to describe the $B_2$-Chebyshev polynomials using a
generalized cosine. Thus consider the change of coordinates
\begin{equation}
    \label{eq:CoordinateChange}
    \begin{pmatrix}
        \theta_1 \\
        \theta_2 
    \end{pmatrix}
    \mapsto
    \begin{pmatrix}
        x_1 \\
        x_2 
    \end{pmatrix}
    =
    \begin{pmatrix}
        \cos(2 \pi \theta_2) \cos(2 \pi ( \theta_1 - \theta_2)) \\
        \cos(\pi \theta_1) \cos( \pi (\theta_1 - 2 \theta_2))
    \end{pmatrix},
\end{equation}
which maps an isoceles right triangle to a surface bounded one, by two
lines and a parabola. One possible choice of the isoceles right triangle
$F$, where this coordinate change is one-to-one, has vertices
$V(F) = \left\{
\begin{pmatrix}
    0 \\
    0
\end{pmatrix},
\begin{pmatrix}
    \frac{1}{2} \\
    \frac{1}{2}
\end{pmatrix},
\begin{pmatrix}
    0 \\
    \frac{1}{2}
\end{pmatrix}
\right\}$, i.e. one has $F = \{ \theta \in \mathbb{R}^2 \; | \;  0 \leq \theta_1
\leq \theta_2 \leq \frac{1}{2} \}$.

The bivariate Chebyshev polynomials of type $B_2$ are
defined using these coordinates as
\begin{equation}
    \label{eq:ChebyshevB2}
    \begin{split}
    T_{k,\ell}(x_1,x_2)
    = & \frac{1}{4} \big(\cos(2 \pi ( k \theta_1 + \ell \theta_2)) \\
    &+ \cos(2 \pi ( (k + \ell) \theta_1 - \ell \theta_2)) \\
    &+ \cos( 2 \pi (k \theta_1 - (2 k + \ell) \theta_2)) \\
    &+ \cos (2 \pi ( (k + \ell) \theta_1 - (2 k + \ell) \theta_2)) 
    \big),         
    \end{split}
\end{equation}
which specializes to
\begin{equation}
    \label{eq:ChebyshevB2OneZero}
    \begin{split}
        T_{n,0}
        &= \cos(2 n \pi \theta_2) \cos(2 n \pi ( \theta_1 - \theta_2)) \\ 
        T_{0,n}
        &= \cos( n \pi \theta_1) \cos( n \pi (\theta_1 - 2 \theta_2)). 
    \end{split}
\end{equation}
Using \eqref{eq:ChebyshevB2OneZero} one can show by elementary
calculations that $T_{n,0}$ and $T_{0,n}$ have $\tfrac{n (n+1)}{2}$
common zeros in $F$. First observe that $\cos( 2 \pi n \theta)$
vanishes for $\theta \in \pm \tfrac{1}{4n} + \tfrac{1}{n} \mathbb{Z}$
and $\cos( \pi n \theta)$ vanishes for
$\theta \in \pm \tfrac{1}{2n} + \tfrac{2}{n} \mathbb{Z}$. Then one has
to choose those $(\theta_1, \theta_2) \in F$. These considerations
result to the common zeros in $F$ being in
$(\theta_1,\theta_2)$-coordinates
\begin{equation}
    \label{eq:B2CommonZerosThetaForm}
    \{ (\tfrac{k}{2n}, \tfrac{j}{4n})
    \; | \; k = 0,\dots,n-1;  j = 1,3,\dots,2n-1;  j \geq 2 k \}.    
\end{equation}
In fact these common zeros are even common zeros for all
$\{T_{k,\ell} \; | \; k+\ell=n\}$. In the sequel we will denote these
common zeros by $\{ \alpha = (\alpha_1, \alpha_2)\}$.

The bivariate Chebyshev polynomials of type $B_2$ are subject to the
following recurrence relations
\begin{equation}
    \label{eq:RecurrenceRelationChebyshevB2}
    \begin{split}
        x_1 \cdot T_{k,\ell} &= \tfrac{1}{4} \left( T_{k+1,\ell} + T_{k-1,\ell}
            + T_{k-1,\ell+2} + T_{k+1,\ell-2} \right), \\
        x_2 \cdot T_{k,\ell} &= \tfrac{1}{4} \left( T_{k,\ell+1} + T_{k,\ell-1}
            + T_{k-1,\ell+1} + T_{k+1,\ell-1} \right).
    \end{split}
\end{equation}
Furthermore they have the decomposition property, i.e.\,
\begin{equation}
    \label{eq:SemigroupPropertyChebyshevB2}
    T_{k,\ell}(T_{n,0}, T_{0,n}) = T_{n \cdot k, n \cdot \ell}
\end{equation}
for all $k,\ell \in \mathbb{N}_0$.

Recall that an ideal $I$ is radical if $f^m \in I$ then $f \in I$. In
the univariate case this corresponds to the polynomial $p(x)$ being
square-free. The number of common zeros and the dimension of
$\mathbb{C}[x_1,x_2] \big/ \ideal{T_{n,0}, T_{0,n}}$ only coincide if
the ideal $\ideal{T_{n,0}, T_{0,n}}$ is radical. In the case of
$B_2$-Chebyshev polynomials this is unfortunately not the case. But
one can always take the radical of an ideal
$\sqrt{I} = \{f \in \algebra{A} \; | \; f^m \in I \text{ for some } m
\in \mathbb{N}\}$.

Now we have all the data we need to define an algebraic signal model.
Consider the polynomial algebra
$\algebra{A} = \mathbb{C}[x_1,x_2] \big/ \sqrt{\ideal{T_{n,0},
    T_{0,n}}}$, the regular module $M = \algebra{A}$, by the common
zeros of $T_{n,0}$ and $T_{0,n}$ of dimension $\tfrac{n(n+1)}{2}$,
with basis $\{T_{k,\ell} \; | \; k + \ell < n\}$, and the
$z$-transform
$\Phi \colon s \mapsto \sum_{k+\ell < n} s_{k,\ell} T_{k,\ell}$. Since
we have determined the common zeros of $T_{n,0}$ and $T_{0,n}$ the
definition of the Fourier transform for this signal model is straight
forward as
\begin{equation}
    \label{eq:FourierTransformTriangle}
    \mathbb{C}[x_1,x_2] \big/ \sqrt{\ideal{T_{n,0}, T_{0,n}}} \cong
    \bigoplus_{\alpha} \mathbb{C}[x_1,x_2] \big/ \ideal{x_1 - \alpha_1, x_2 -
      \alpha_2}.
\end{equation}
It can be realized via choice of $(1)$ as a basis in
each $\mathbb{C}[x_1,x_2] \big/ \ideal{x_1 - \alpha_1, x_2 - \alpha_2}$ by the
matrix
\begin{equation}
    \label{eq:FourierTransformMatrixTriangle}
    \Fourier_n = (T_{k,\ell}(\alpha))_{k+\ell < n, \alpha}.
\end{equation}

\begin{figure}
    \centering
    \includegraphics[width=0.33\textwidth]{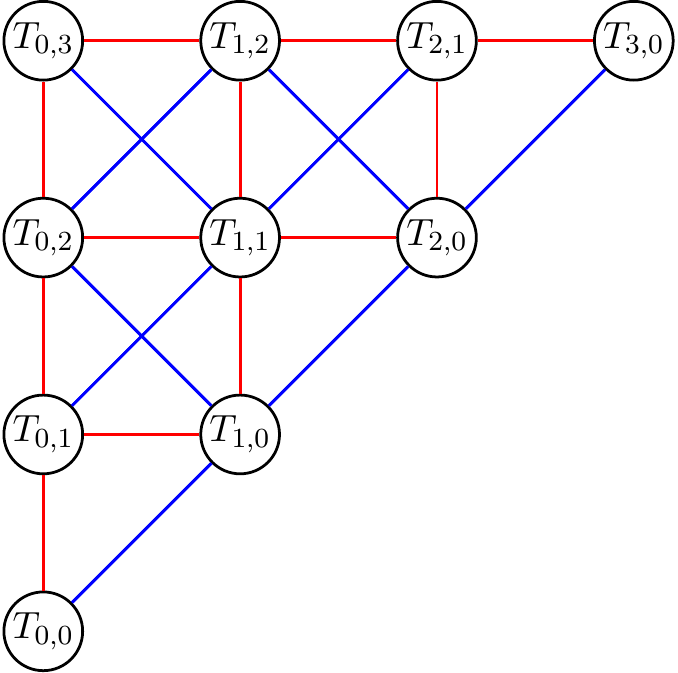}
    \caption{Visualization of the triangle signal model. The
      $x_1$-shifts are blue colored and the $x_2$-shifts are red colored.} 
    \label{fig:TriangleSignalModel}
\end{figure}

The recurrence relations~\eqref{eq:RecurrenceRelationChebyshevB2}
define the structure of the visualization graph in
Fig.~\ref{fig:TriangleSignalModel}. 

\section{Christoffel-Darboux and the inverse transform}
\label{sec:ChristoffelDarbouxForInverse}%

It would be quite nice if one has a unitary or orthogonal version of a
discrete transform, since then one only needs to find and implement
one algorithm for the computation of the transform and the computation
of the inverse transform. It is well-known that in the case of 1D
signal transforms the existence of unitary versions is connected to
the basis polynomials of the signal module being orthogonal
polynomials and relys on the Christoffel-Darboux formula for
univariate orthogonal polynomials~\cite{Yemini.Pearl:1979a}.

In the multivariate setting there is another constraint - the
vanishing of all basis polynomials of the same degree. To see this, we
recall the multivariate Christoffel-Darboux formula from
\cite{Xu:1993a}. Denote by
$\mathbb{T}_k = (T_{0,k}, T_{1,k-1}, \dots, T_{k,0})^\top$ the vector
of bivariate Chebyshev polynomials of degree $k$. The vector
$\mathbb{T}_k$ is of length $k+1$. As the Chebyshev polynomials are
orthogonal polynomials they satisfy three-term recurrence relations
\begin{equation}
    \label{eq:ThreeTermRecurrenceRelation}
    x_i \mathbb{T}_k = A_{k,i} \mathbb{T}_{k+1} + B_{k,i} \mathbb{T}_k
    + C_{k,i} \mathbb{T}_{k-1},
\end{equation}
were the matrices $A_{k,i}, B_{k,i},$ and $C_{k,i}$ can be deduced
from~\eqref{eq:RecurrenceRelationChebyshevB2}. For example from the
$x_1$-shift one gets the matrices
\begin{equation}
    \label{eq:ThreeTermRecurrenceRelationMatrices}
    \begin{split}
            A_{k,1} &=
            \begin{bsmallmatrix}
                0 & 1/2 & 0 & \dots & & 0 \\
                1/4 & 0 & 1/4 & 0 & \dots & 0 \\
                0 & \ddots & \ddots & \ddots &  & \vdots \\
                \dots & 0 & 1/4 & 0 & 1/4 & 0 \\
                0 & \dots & 0 & 1/2 & 0 & 1/4 
            \end{bsmallmatrix}, \\
            B_{k,1} &=
            \begin{bsmallmatrix}
                0 & \dots &  & 0 \\
                \vdots & \ddots & &\vdots \\
                0 & \dots & 0 & 0 \\
                0 & \dots & 1/4 & 0 \\
                0 & \dots & 0 & 0
            \end{bsmallmatrix},
            C_{k,1} =
            \begin{bsmallmatrix}
                0 & 1/2 & 0 & \dots & 0 \\
                1/4 & 0 & 1/4 & \ddots & \vdots  \\
                0 & \ddots & \ddots & \ddots &   0 \\
                & & 1/4 & 0 & 1/4 \\
                \vdots & & 0 & 1/4 & 0 \\
                0 & \dots & & 0 & 1/4 
            \end{bsmallmatrix},
    \end{split}
\end{equation}
with special case $B_{1,1} =
\begin{bsmallmatrix}
    1/2 & 0 \\
    0 & 0 
\end{bsmallmatrix}$.

From the three-term recurrence relation one can deduce a multivariate
Christoffel-Darboux formula~\cite{Xu:1993a} 
\begin{equation}
    \label{eq:ChristoffelDarbouxFormula}
    \begin{split}
    &\sum_{k=0}^{n-1} \mathbb{T}_k^{\top}(x) H_k^{-1} \mathbb{T}_k(y) \\
    &= 
    \begin{cases}
        \begin{split}
            &(x_i - y_i)^{-1} \cdot \\
            &\big((A_{n-1,i} \mathbb{T}_{n}(x))^\top
            H_{n-1}^{-1}
            \mathbb{T}_{n-1}(y)  \\
            &- \mathbb{T}_{n-1}^\top(x) H_{n-1}^{-1} A_{n-1,i}
            \mathbb{T}_{n}(y) \big)
        \end{split}
        &\text{if } x_i \not= y_i \\
        \begin{split}
            &\mathbb{T}_{n-1}^\top(x) H_{n-1}^{-1} A_{n-1,i}
            \tfrac{\partial}{\partial x_i} \mathbb{T}_{n}(x) \\
            &- (A_{n-1,i} \mathbb{T}_{n}(x))^\top H_{n-1}^{-1}
            \tfrac{\partial}{\partial x_1} \mathbb{T}_{n-1}(x)
        \end{split}
        &\text{if } x_i = y_i,
    \end{cases}        
    \end{split}    
\end{equation}
with matrices $H_0 = \tfrac{1}{2}$ and
$H_k = \mathsf{diag}(\tfrac{1}{8}, \tfrac{1}{16}, \dots,
\tfrac{1}{16}, \tfrac{1}{8})$.

Now one can observe that the entries of the matrix $\Fourier_n^\top
\cdot \Fourier_n$ are of the form $\sum_{k=0}^{n-1}
\mathbb{T}_k^\top(\alpha) \mathbb{T}_k(\beta)$ for $\alpha,\beta$
common zeros of $T_{n,0},T_{0,n}$. Thus consider the matrix
$H_n^{\oplus} = \bigoplus_{k=0}^{n-1} H_k^{-1}$. Since the common zeros of
$T_{n,0}$ and $T_{0,n}$ are in fact common zeros of all entries of
$\mathbb{T}_n$ one obtains the diagonal matrix
\begin{equation}
    \label{eq:FourierTimesHWeightsGivesDiagonal}
    \Fourier_n^\top \cdot H_n^{\oplus} \cdot \Fourier_n
    = \mathsf{diag}\left(\mathbb{T}_{n-1}^\top(\alpha) H_{n-1}^{-1} A_{n-1,1}
    \tfrac{\partial}{\partial x_1} \mathbb{T}_{n}(\alpha)\right). 
\end{equation}
Since the diagonal entries of
\eqref{eq:FourierTimesHWeightsGivesDiagonal} do not vanish we can
invert them. Denote the diagonal matrix with inverted entries by
\begin{equation}
    \label{eq:DiagonalMatrixChristoffelEntries}
    D_n = \mathsf{diag}\left( \left(\mathbb{T}_{n-1}^\top(x)
            H_{n-1}^{-1} A_{n-1,1} \tfrac{\partial}{\partial x_1}
            \mathbb{T}_{n}(x)\right)^{-1}\right). 
\end{equation}
One obtains
\begin{equation}
    \label{eq:InverseB2TrafoByChristoffelDarboux}
    \Fourier_n^{-1} = D_n \Fourier_n^\top H_n^{\oplus},
\end{equation}
and in turn an orthogonal version of the triangle transform
\begin{equation}
    \label{eq:OrthogonalB2TrafoByChristoffelDarboux}
    \Fourier_n^{\mathsf{orth}} = \sqrt{H_n^{\oplus}} \Fourier_n
    \sqrt{D_n}. 
\end{equation}

\section{Conclusions and future work}
\label{sec:Conclusion}%

We developed a novel cosine transform on triangles and derived an
orthogonal version of it. This showed the significance of multivariate
Chebyshev polynomials and the multivariate Christoffel-Darboux formula
for the derivation of orthogonal transforms.

To make the usage of these transforms applicable in real world
applications a fast algorithm is essential. Since the multivariate
Chebyshev polynomials obey a decomposition
property~\eqref{eq:SemigroupPropertyChebyshevB2}, a fast algorithm
exists. The implementation of this fast algorithm is currently under
consideration and will be subject of upcoming work.


\vfill\pagebreak

\bibliographystyle{IEEEbib}
\bibliography{ICASSP}

\end{document}